\documentclass[10pt,a4paper]{article}
\usepackage{titlesec}
\usepackage{enumitem}
\usepackage{url}
\usepackage[utf8]{inputenc}
\titleformat*{\section}{\centering\bfseries\large}
\titleformat*{\subsection}{\centering\bfseries}
\usepackage{amsmath}
\usepackage[arrow, matrix, curve]{xy}
\usepackage{amssymb}
\usepackage{a4wide}
\usepackage{tikz}
\usepackage{amsthm}
\usepackage{latexsym}
\newtheorem{thm}{Theorem}[section]
\newtheorem{lem}[thm]{Lemma}
\newtheorem{cor}[thm]{Corollary}

\theoremstyle{definition}

\theoremstyle{remark}
\newtheorem{remark}[thm]{Remark}
\newtheorem{example}[thm]{Example}

\DeclareMathOperator{\PGl}{PGL}
\DeclareMathOperator{\PSl}{PSL}

\DeclareMathOperator{\nrd}{Nrd}

\newcommand{\CC}{\mathbb C}

\newcommand{\PP}{\mathbb P}
\newcommand{\QQ}{\mathbb Q}
\newcommand{\RR}{\mathbb R}

\newcommand{\HH}{\mathbb H}

\newcommand{\pP}{\mathfrak p}

\newcommand{\oo}{\mathfrak o}
\newcommand{\OO}{\mathfrak O}

\newcommand{\cC}{\mathcal C}
\begin{document}
\title{Lower bounds for covolumes of arithmetic lattices in $\PSl_2(\RR)^n$}
\author{Amir D\v{z}ambi\'c\footnote{Institut f\"ur Mathematik, Goethe-Universit\"at Frankfurt a.M., Postfach 111932, Fach 187, 60054 Frankfurt a.M. Email: dzambic@math.uni-frankfurt.de}}
\maketitle
\begin{abstract}
\noindent We study the covolumes of arithmetic lattices in $\PSl_2(\RR)^n$ for $n\geq 2$ and identify uniform and non-uniform irreducible lattices of minimal covolume. More precisely, let $\mu$ be the Euler-Poincar\'e measure on $\PSl_2(\RR)^n$ and $\chi=\mu/2^n$. We show that the Hilbert modular group $\PSl_2(\oo_{k_{49}})\subset \PSl_2(\RR)^3$, with $k_{49}$ the totally real cubic field of discriminant $49$ has the minimal covolume with respect to $\chi$ among all irreducible lattices in $\PSl_2(\RR)^n$ for $n\geq 2$ and is unique such lattice up to conjugation. The cocompact lattice of minimal covolume with respect to $\chi$ is the normalizer $\Delta_{k_{725}}^u$ of the norm-1 group of a maximal order in the quaternion algebra over the unique totally real quartic field with discriminant $725$ ramified exactly at two infinite places, which is a lattice in $\PSl_2(\RR)^2$. There is exactly one more lattice in $\PSl_2(\RR)^2$ and exactly one in $\PSl_2(\RR)^4$ with the same covolume as $\Delta_{k_{725}}^u$, which are the Hilbert modular groups corresponding to $\QQ(\sqrt{5})$ and $k_{725}$. The two lattices $\Delta_{k_{725}}^u$ and $\PSl_2(\oo_{\QQ(\sqrt{5})})$ have the smallest covolume with respect to the Euler-Poincar\'e measure among all arithmetic lattices in $G_n$ for all $n\geq 2$. These results are in analogy with Siegel's theorem on the unique minimal covolume (uniform and non-uniform) Fuchsian groups and its generalizations to various higher dimensional hyperbolic spaces due to Belolipetsky, Belolipetsky-Emery, Stover and Emery-Stover. 

\end{abstract}
\section{Introduction and the statement of the main result}
Let $G$ be a semisimple real Lie group with trivial center and without compact factors, and let $\mu$ be a Haar measure on $G$. A lattice $\Gamma\subset G$ is a discrete subgroup of $G$ such that the covolume $\mu(G/\Gamma)$ is finite. A lattice is called uniform or cocompact if $G/\Gamma$ is compact. By a result of Ka\v{z}dan-Margulis \cite{KazdanMargulis}, there exist a constant $C=C_{G,\mu}>0$ such that $\mu(G/\Gamma)\geq C$ for any lattice $\Gamma$ in $G$. It is an open and interesting problem to find explicit non-trivial lower bounds and minimal values for $\mu(G/\Gamma)$, where $\Gamma$ runs through the set of all lattices in $G$, for a given Lie group $G$. 
The first result in this direction is the theorem of Siegel who proved that for $G=\PSl_2(\RR)$ and $\mu$ the Euler-Poincar\'e measure we have $\mu(G/\Gamma)\geq 1/42$ and $\mu(G/\Gamma)=1/42$ if and only if $\Gamma=\Delta(2,3,7)$ is the Fuchsian triangle group of signature $(2,3,7)$, which turns out to be an arithmetic group. The analogous  result for $G=\PSl_2(\CC)$ has been first established only for arithmetic lattices \cite{ChinburgFriedman86} and only after a long work has been settled without the assumption on arithmeticity \cite{GehringMartin1},\cite{GehringMartin2}. Again, there is a minimum $min_{\Gamma\subset G}\mu(G/\Gamma)$, which is attained by a unique lattice. Similar results have been established in the case of the group $PO(n,1)$ and $PU(n,1)$ see \cite{AdeboyeWei12}, \cite{AdeboyeWei14} and references therein. Note that for general $n$ the results are weaker than those established for $G=\PSl_2(\RR)$ and $\PSl_2(\CC)$ as we do not know if the lower bounds are indeed attained by some lattice. Such informations are available only for arithmetic lattices, see \cite{BelolipetskyEmery}, \cite{EmeryStover14}, \cite{Stover11}, \cite{Stover14}.\\

In the present note we will discuss the question of finding the lower bounds for covolumes of lattices in $G_n:=\PSl_2(\RR)^n$. To be more precise, let $\mu$ denote the Euler-Poincar\'e measure on $G_n$, see \cite{Serre71}, which is the unique Haar measure on $G_n$ such that for a torsion-free and uniform lattice $\Gamma\subset G_n$ the covolume $\mu(\Gamma):=\mu(G_n/\Gamma)=(-1)^n e(\Gamma\backslash \HH^{(n)})$ is up to the sign the Euler number of the corresponding locally symmetric space $\Gamma\backslash\HH^{(n)}$, where $\HH^{(n)}=\HH\times\ldots\times \HH$ is an $n$-fold product of complex upper half planes $\HH$. Note that $\mu$ can be realized as the volume of a fundamental domain of $\Gamma$ on $\HH^{(n)}$ with respect to the volume form $(-1)^nd\nu/(2\pi)^n$ where $d\nu=\prod_{1\leq i\leq n} dx_i\wedge dy_i/y_i^2$ is the standard hyperbolic volume form on $\HH^{(n)}$. For computational reasons we will initially use another normalization of the Haar measure and consider the $G_n$ invariant measure $\chi=\mu/2^n$, which for torsion-free and uniform $\Gamma\subset G_n$ has also a geometric meaning as $\chi(\Gamma)=(-1)^n\chi(\Gamma\backslash \HH^{(n)})$, where $\chi(\Gamma\backslash \HH^{(n)})$ is the holomorphic Euler characteristic (arithmetic genus) of the corresponding locally symmetric space, which follows from the work of Matsushima-Shimura \cite{MatsushimaShimura63}. 

For $n\geq 2$ one has to distinguish between two types of lattices in $G_n$; A lattice $\Gamma\subset G_n$ is called irreducible if the image of $\Gamma$ under the canonical projection $G_n\longrightarrow G_{n'}$ for any $1\leq n'<n$ is not a discrete subgroup in $G_{n'}$. Otherwise, $\Gamma$ contains with finite index a product $\Gamma_1\times \Gamma_2$ of two lattices $\Gamma \subset G_{n'}$ and $\Gamma_2\subset G_{n-n'}$. In this case $\Gamma$ is called reducible. Since the rank of $G_n$ is $n$, the famous arithmeticity theorem of Margulis implies that for $n\geq 2$ every irreducible lattice $\Gamma$ in $G_n$ is arithmetic. This means that $\Gamma$ is commensurable with the group of norm-1-elements of an order in a quaternion algebra $B$ over a totally real number field $k$ of degree $m$ ramified at $m-n$ infinite places of $k$ (see Section \ref{prelim} for the precise definition).    
 
Our main result is the following 
\begin{thm}
\label{mainthm}
Let $n\geq 2$ be an integer and $G_n=\PSl_2(\RR)^n$. Let $\mu$ be the Euler-Poincar\'e measure on $G_n$ and let $\chi=\mu/2^n$ be the Haar measure on $G_n$ such that $\chi(\Gamma):=\chi(G_n/\Gamma)$ equals  $(-1)^n\chi(X_{\Gamma})$ where $\chi(X_{\Gamma})$ is the holomorphic Euler characteristic of the corresponding variety $X_{\Gamma}=\Gamma\HH^{(n)}$ for $\Gamma$ torsion-free and uniform. Then: 
\begin{enumerate}
\item For any $n\geq 2$ and any arithmetic lattice $\Gamma\subset G_n$ we have $\chi(\Gamma)\geq 1/84$. The minimal value $1/84$ is attained exactly by one lattice, which is
the (non-uniform) Hilbert modular group $\Delta_{k_{49}}^{nu}:=\PSl_2(\oo_{k_{49}})\subset G_3$. 
\item For $\Gamma\subset G_n$ uniform we have $\chi(\Gamma)\geq 1/60$. This minimal value is attained by the unique (up to conjugation) uniform lattice $\Delta_{k_{725}}^{u}\subset G_2$, the maximal lattice defined by the quaternion algebra $B$ over the totally real quartic field of discriminant $725$ ramified at exactly two infinite places. The minimal value $1/60$ is also attained by exactly two non-uniform lattices, which are the Hilbert modular group $\Delta_{k_5}^{nu}:=\PSl_2(\oo_{\QQ(\sqrt{5})})\subset G_2$ of the real quadratic field $\QQ(\sqrt{5})$ and the Hilbert modular group $\Delta_{k_{725}}^{nu}=\PSl_2(\oo_{k_{725}})\subset G_{4}$.  
\item The minimal covolume $\mu(\Gamma)$ of an arithmetic lattice $\Gamma\subset G_n$ with respect to $\mu$ is $1/15$. This value is attained exactly by the two lattices $\Delta_{k_{725}}^u$ and $\Delta_{k_5}^{nu}$.  
\end{enumerate}
\end{thm} 
For the precise definition of the lattice $\Delta_{k_{725}}^u$ see Sections \ref{prelim} and \ref{examples}.\\
     
Note the remarkable property that the minimal covolume in $G_2=\PSl_2(\RR)^2$ is attained by exactly two lattices, one uniform and one non-uniform. This example seems to contradict the tendency that the minimal covolume is attained by a unique lattice, or lattices which are either uniform or non-uniform, known from other examples. Also note that the value $1/84$ is attained by the known Fuchsian group of minimal covolume, which is arithmetic and defined by a quaternion algebra over the totally real cubic field of discriminant $49$. We should note that the existence of a universal lower bound for the covolume in $G_n$, which holds for all $n$ is restricted only to irreducible lattices. Namely, as for a direct product $\Gamma_1\times \Gamma_2$ of lattices $\Gamma_1\subset G_{n_1}$, $\Gamma_{n_2}$ the relation $\chi(\Gamma_1\times \Gamma_2)=\chi(\Gamma_1)\chi(\Gamma_2)$ holds, taking any lattice $\Delta\subset G_1$ with $\chi(\Delta)<1$ we see that the covolume $\chi(\Delta^n)$ of the reducible lattice $\Delta^n\subset G_n$ tends to zero as $n$ goes to infinity.\\
 
Let us also mention an immediate consequence of the second part of the Theorem \ref{mainthm}, which was one motivation to study the lattices of small covolumes:\\
Recall that for any uniform and torsion-free lattice $\Gamma\subset G_n$ the orbit space $X_{\Gamma}=\Gamma\backslash \HH^{(n)}$ has a structure of a smooth $n$-dimensional projective variety of general type. In the case where $\Gamma$ is arithmetic we call such a variety a \emph{quaternionic Shimura variety}. As it is well-known, the automorphism group of such a variety is finite. More precisely, the automorphism group of $X_{\Gamma}$ is isomorphic to the factor group $N\Gamma/\Gamma$, where $N\Gamma$ is the normalizer of $\Gamma$ in $G_n$, which is again arithmetic. Now, for arbitrary finite index inclusion $\Gamma\subset \Delta$ of lattices in $G_n$ we have the equality $\chi(\Gamma)=[\Delta:\Gamma]\chi(\Delta)$, which implies that $|Aut(X_{\Gamma})|\leq \chi(\Gamma)/C^a_{G_n,\chi}$, where $C^a_{G_n,\chi}$ is any lower bound for the covolume of a uniform arithmetic lattice in $G_n$. By the main Theorem $C^a_{G_n,\chi}=1/60$ and we get:
\begin{cor}
For a smooth compact quaternionic Shimura variety $X_{\Gamma}=\Gamma\backslash \HH^{(n)}$ the order of the automorphism group is bounded by $60\chi(X_{\Gamma})$, where $\chi(X_{\Gamma})$ denotes the arithmetic genus of $X_{\Gamma}$. 
\end{cor}
 The quaternionic Shimura varieties with $\chi(X_{\Gamma})=1$ are of particular interest as they give examples of smooth projective varieties with the same Betti numbers as $(\PP^1)^n$, see \cite{Dzambic14}, \cite{Dzambic14a}. The corollary implies that the number of automorphisms of such a variety is bounded by $60$.     

\section{Maximal arithmetic lattices in $\PSl_2(\RR)^n$ and volumes}

\subsection{Preliminaries}
\label{prelim}
Let $k$ be a totally real field of degree $m=[k:\QQ]$, $\oo_k$ its ring of integers, and $U_k$ the group of units in $k$. An element $x\in k$ is totally positive if $\sigma(x)$ is a positive real number for every embedding $\sigma:k\longrightarrow \RR$. We put $k_+=\{x\in k\mid x\ \text{is totally positive}\}$ and $\oo_{k,+}=\oo_k\cap k_+$, $U_{k}^+=U_k\cap k_+$. Let $V_k$ denote the set of all places of $k$, that is, the set of equivalence classes of absolute values on $k$. Then $V_k=V^{\infty}_{k}\cup V_k^f$, where $V_k^{\infty}$ denotes the set of infinite places and $V_k^f$ the set of all finite places of $k$. Each $v_{j,\infty}\in V_k^{\infty}$ corresponds to a field embedding $k\longrightarrow \RR$ whereas a finite place $v\in V_k^f$ corresponds to a non-trivial prime ideal of $\oo_k$.  Let $B$ be a quaternion algebra over $k$. For each $v\in V_k$ let $B_{v}=B\otimes_{k}k_v$ be the quaternion algebra over the completion $k_v$ of $k$ with respect to $v$ arising from $B$ by scalar extension. We say that $B$ is \emph{ramified} at $v$ if $B_v$ is a division algebra and \emph{unramified} otherwise (i.e. $B_v\cong M_2(k_v)$). Put $R(B)=\{v\in V_k\mid B\ \text{is ramified at}\ v\}=R_{\infty}(B)\cup R_f(B)$ with $R_{\infty}(B)=R(B)\cap V_k^{\infty}$ and $R_{f}(B)=R(B)\cap V_k^f$. From general theory of central simple algebras we know that $R(B)$ is a finite set of even order and moreover we know that $B$ is isomorphic to $M_2(k)$ if and only if $R(B)=\emptyset$. Otherwise, that is, if $|R(B)|\neq 0$, the algebra $B$ is a division algebra. Let $B^+=\{x\in B\mid \nrd(x)\in k_+\}$ and for a maximal order $\OO$ in $B$ let 
$$
\OO^+=\OO\cap B^+=\{x\in \OO\mid \nrd(x)\in \oo_{k,+}\}
$$
where $\nrd(\ )$ denotes the reduced norm on $B$. Assume that $|R_{\infty}(B)|=n$. Then $n\leq m$ and the diagonal embedding $B\hookrightarrow M_2(\RR)^n$ induces an embedding $\Gamma_{\OO}^+:=\OO^+/U_k\hookrightarrow \PSl_2(\RR)^n=:G_n$. The image of this embedding is an irreducible lattice in $G_n$. More precisely, for any projection $\pi_s:G_n\longrightarrow G_s$ onto a non-trivial partial product of $G_n$ the image $\pi_s(\Gamma)$ is a dense (in particular non-discrete) subgroup. 
We define $\cC(k,B)$ as the set of all lattices in $G_n$ which are commensurable with $\Gamma_{\OO}^+$ where $\OO$ is some maximal order in the quaternion algebra $B$ over $k$ with the above properties. We say that a lattice $\Gamma\in G_n$ is \textbf{arithmetic} if there is a field $k$ and quaternion algebra $B$ as above such that $\Gamma\in \cC(k,B)$. By the Godement's criterion, $\Gamma\in \cC(k,B)$ is uniform (that is, $G_n/\Gamma$ is compact) if and only if $B$ is a division algebra, that is, $B\ncong M_2(k)$. Hence, $\Gamma$ is non-uniform if and only if $\Gamma$ belongs to $\cC(k,M_2(k))$. By the above discussion in the non-uniform case we necessarily have $m=n$. \\
As a locally compact group $G_n$ possesses a Haar measure (unique up to a multiplicative constant), and moreover $G_n$ is unimodular, that is, the left invariant Haar measure is also right invariant. As already mentioned, in the following we will work with two particular normalizations of the Haar measure, which are denoted by $\mu$ and $\chi$: If $\Gamma$ is a torsion-free and uniform lattice in $G_n$ then the value $\mu(\Gamma):=\mu(G_n/\Gamma)$ equals $(-1)^ne(X_{\Gamma})$, with $e(X_{\Gamma})$ the topological Euler characteristic of $X_{\Gamma}=\Gamma\backslash \HH^{(n)}$. The measure $\chi$ is defined as $\chi=\mu/2^n$ and equals $(-1)^n\chi(X_{\Gamma})$, where $\chi(X_{\Gamma})$ is the holomorphic Euler characteristic of $X_{\Gamma}$. We have already used this normalization in \cite{Dzambic14}. Let $\Gamma\in \cC(k,B)$ be arithmetic. Then, by definition, $\Gamma$ is commensurable with a group $\Gamma_{\OO}^1$ associated with a maximal order $\OO$ in $B$ and we define the generalized index $(\Gamma_{\OO}^1:\Gamma)$ as 
$$
(\Gamma_{\OO}^1:\Gamma)=\frac{[\Gamma_{\OO}^1:\Gamma\cap \Gamma_{\OO}^1]}{[\Gamma: \Gamma\cap \Gamma_{\OO}^1]}.
$$  
The generalized index is a rational number and equals to the usual index if $\Gamma<\Gamma_{\OO}^1$.
\begin{thm}(Vign\'eras \cite{Vigneras76})
For any $\Gamma\in \cC(k,B)$ we have 
\begin{equation}
\label{voleins}
\chi(\Gamma)=(\Gamma_{\OO}^1:\Gamma)(-1)^{m+n}2^{n-m-1}\zeta_k(-1)\prod_{v\in R_f(B)}(Nv-1)
\end{equation}
\end{thm} 
By the multiplicativity property of the volume $\chi$, the minimum $\min_{\Gamma\in \cC(k,B)}\{\chi(\Gamma)\}$ will be attained on a maximal group $\Gamma\in \cC(k,B)$. Therefore we need more precise informations on maximal elements in $\cC(k,B)$. Such results on maximal elements of $\cC(k,B)$, where $B$ is a quaternion algebra over a general number field and their covolumes are obtained by A. Borel \cite{Borel81} (see also \cite{Dzambic14} and \cite{Dzambic14a} for the case of lattices in $G_n$). According to Borel we have the following result on the minimal value of $\chi$ on $\cC(k,B)$. 

\begin{thm}(Borel \cite{Borel81})
\label{basic}
\begin{enumerate}
\item Let $B$ be as above and $\OO$ a maximal order in $B$. Let 
$$
N\Gamma_{\OO}^+=\{x\in B^+\mid x\OO x^{-1}=\OO\}/k^{\ast},
$$
be the normalizer of a maximal order $\OO$ in $B^+$. Then 
$$\min_{\Gamma\in \cC(k,B)} \chi(\Gamma)=\chi(N\Gamma_{\OO}^+)$$ 
independently of the choice of $\OO$. 
\item Let $k_B$ be the maximal abelian extension of $k$ unramified at all the finite places of $k$, such that the Galois group $Gal(k_B/k)$ is elementary $2$-abelian and with the property that all the finite places $v\in R_f(B)$ are totally decomposed in $k_B$. Further, let $k'_B$ have the same properties as $k_B$ but with the restriction to be unramified at all the places of $k$. Particularly, $k'_B$ is totally real and 
$$
[k_B:k]=\frac{2^m[k'_B:k]}{[U_k:U_{k}^+]}.
$$
Let $t$ be the number of places of $k$ lying over $2$. Define 
\begin{align}
\label{gkb}
g(k,B)&=\frac{d_k^{3/2}\zeta_k(2)}{2^{2m+t-1}\pi^{2m}[k_B:k]}=\frac{(-1)^m\zeta_k(-1)}{2^{m+t-1}[k_B:k]}\\
&=\frac{d_k^{3/2}\zeta_k(2)[U_k:U_{k}^+]}{2^{3m-1+t}\pi^{2m}[k'_B:k]}
=\frac{(-1)^m\zeta_k(-1)[U_k:U_{k}^+]}{2^{2m+t-1}[k'_B:k]},
\end{align}
and put
$$
\mu(k,B)=2^{t-t'}g(k,B).
$$
Then we have the following exact formula for $\chi(N\Gamma_{\OO}^+)$
\begin{equation}
\label{volmax}
\chi(N\Gamma_{\OO}^+)=\mu(k,B)\prod_{\substack{v\in R_f(B)\\  Nv\neq 2}} \frac{Nv-1}{2}
\end{equation}
where $t'=|\{v\in R_f(B)\mid Nv=2\}|$.
\end{enumerate}
\end{thm}
\subsection{Examples of maximal lattices in $\PSl_2(\RR)^n$ of small volume}
\label{examples}
The above theory particularly applies to the case $B=M_2(k)$ of the matrix algebra over a degree $m$ totally real field $k$, which is characterized by the property $R(B)=\emptyset$. Taking $\OO=M_2(\oo_k)$ which is clearly a maximal order in $M_2(k)$, the arithmetic group $\Gamma_{\OO}^1=\Gamma_{k}=\PSl_2(\oo_k) \subset G_m$ is the well-known \emph{Hilbert modular group}. The group $\Gamma_{\OO}^+=\PGl_2^+(\oo_k)$ is the so-called \emph{extended Hilbert modular group} whereas $N\Gamma_{\OO}^+$ ist the so-called \emph{Hurwitz-Maass extension} of $\Gamma_{k}$. It is not difficult to see that in this case $[\Gamma_{\OO}^+:\Gamma_{\OO}^1]=[U_{k}^+,U_k^2]$ and $[N\Gamma_{\OO}^+:\Gamma_{\OO}^+]=[Cl_k:Cl_k^2]$, where $Cl_k$ is the class group of $k$.

\begin{example}
Let $k_5=\QQ(\sqrt{5})$. Then, as the fundamental unit $\epsilon_{k_5}=\frac{-1+\sqrt{5}}{2}$ is not totally positive, $U_{k_5}^+=U_{k_5}^2$. Moreover, the ideal class number of $k_5$ is one (note that $k'_B$ is a subfield of the Hilbert class field). It follows that the Hilbert modular group $\Delta_{k_5}^{nu}:=\PSl_2(\oo_{k_5})$ is maximal, that is, coincides with its Hurwitz-Maass extension. Using formula (\ref{voleins}) (or formula (\ref{volmax})) we find that $\chi(\Delta_{k_5}^{nu})=\zeta_{k_5}(-1)/2$. On the other hand, we can use the Siegel-Klingen formula \cite{Siegel1969} for the value $\zeta_{k_5}(-1)$ by which we have $\chi(\Delta_{k_5}^{nu})=1/60$.
\end{example}

\begin{example}
Now, let $k_{49}$ be the unique totally real cubic number field of discriminant $49$. The defining polynomial of $k_{49}$ is $X^3-2X^2-X+1$ and using the computer algebra system SAGE we find that the fundamental units of $k_{49}$ are not totally positive. It follows that $U_{k_{49}}^+=U_{k_{49}}^2$. Also in this case SAGE computes the ideal class number of $k_{49}$, which is one. It follows that the Hilbert modular group $\Delta_{k_{49}}^{nu}=\PSl_2(\oo_{k_{49}})$ is maximal. Since prime $2$ is inert in $k_{49}$, again formula (\ref{volmax}) implies that $\chi(\Delta_{k_{49}}^{nu})=-\zeta_{k_{49}}(-1)/4$. After identifying the value $\zeta_{k_{49}}(-1)=1/21$, we have $\chi(\Delta_{k_{49}}^{nu})=1/84$.
\end{example}

\begin{example}
Consider now the unique totally real quartic field $k_{725}$ of discriminant $d_{k_{725}}=725$. The field $k_{725}$ is defined by the irreducible polynomial $X^4-X^3-3X^2+X+1$. Using SAGE again, we compute the fundamental units of $k$ and conclude that $U_k^+=U_k^2$ in this case. Also, using SAGE for instance we know that the class number of $k_{725}$ is one. 
Computing the value $\zeta_{k_{725}}(-1)=2/15$ using SAGE and the already mentioned method of Siegel-Klingen, we can deduce the value $\chi(\Delta_{k_{725}}^{nu})=\zeta_{k_{725}}(-1)/8=1/60$ for the Hilbert modular group $\Delta_{k_{725}}^{nu}:=\PSl_2(\oo_{k_{725}})$.\\
We consider now the quaternion algebra $B$ over $k_{725}$ which is unramified exactly at two infinite places of $k_{725}$. Let us consider an arbitrary maximal order $\OO$ in $B$. Since $R_f(B)=\emptyset$, the formula (\ref{volmax}) gives us for $\Delta_{k_{725}}^{u}=N\Gamma_{\OO}^+$ the value $\chi(\Delta_{k_{725}}^{u})=\zeta_{k_{725}}(-1)/8$.  We obtain $\chi(\Delta_{k_{725}}^{u})=1/60$.   
\end{example}

\section{Proof of the main Theorem}

In this section we consider an arbitrary arithmetic lattice $\Gamma\subset G_n$ for $n\geq 2$ with $\chi(\Gamma)\leq 1/60$. Since $\Gamma$ is arithmetic, it belongs to a commensurability class $\cC(k,B)$ where $k$ is a totally real number field of degree $m\geq n$ and $B$ is a quaternion algebra over $k$ ramified at $m-n$ infinite places of $k$. By Theorem \ref{basic}, we know that $\chi(\Gamma)\geq \chi(N\Gamma_{\OO}^+)$ with some maximal $N\Gamma_{\OO}^+$, hence we know that 
$$
\frac{1}{60}\geq \mu(k,B)\prod_{\substack{v\in R_f(B)\\  Nv\neq 2}} \frac{Nv-1}{2}.
$$ 
Since $t\geq t'$ and $\frac{Nv-1}{2}>1$ for all $v\in V_k$ such that $Nv\neq 2$, we obtain the inequality 
\begin{equation}
\label{1/60}
g(k,B)\leq \mu(k,B)\leq \frac{1}{60}.
\end{equation}
We will now follow the strategy of Chinburg-Friedman \cite{ChinburgFriedman86} to show in several steps that the inequality (\ref{1/60}) implies that the degree of $k$ is less or equal $5$, that the degrees $[k'_B:k]$ and $[k_B:k]$ have to be one. Using this information, in the last step we will analyze all possibilities for pairs $(k,B)$ with $m\leq 5$ individually to finally get the stated result.  

\subsection{Basic inequalities}
\label{inequalities}
Let us recall the lower bounds for $g(k,B)$ and $\mu(k,B)$ due to Chinburg-Friedman.

\begin{thm}(\cite[Proposition 3.1]{ChinburgFriedman86})
Let $B$ and $k$ be as before and let $K$ be some finite Galois extension of $k$ unramified at all the finite places of $k$. Put $[K:\QQ]=m_K=r_1(K)+2r_2(K)$, with $r_1(K)$ the number of real places of $K$. Then for every $1<s\leq 2$ and $y>0$ 
\begin{equation}
\label{chfr1}
\log g(k,B) > \log(0.0124[U_k:U_{k}^+])-\log(s(s-1))+m(0.738-\log(\Gamma(s/2)))+mT(s,y)
\end{equation}
where $T(s,y)=T(s,y,\mathcal S,p_0)$ is the function defined by 

\begin{align*}
T(s,y)&= \frac{3\gamma}{2}+\log(2\sqrt{\pi})+\frac{3r_1(K)}{2m_K}-\log(\pi)\\
&-s\left( \frac{\gamma+\log(2)}{2}+\frac{r_1(K)}{2m_K} +\frac{\log(2)}{2}\right)\\
&-\frac{3-s}{2}\int_{0}^{\infty} (1-a(x\sqrt{y}))\left\{  \frac{1}{\sinh(x)}+\frac{r_1(K)}{2m_K\cosh^2(x/2)}\right\}dx\\
&-\frac{(3-s)6\pi}{5m_K\sqrt{y}}\\
&+R(s,p_0)\\
&+\frac{1}{m_K}\sum_{p\in \mathcal S}\sum_{\mathfrak p\mid p} q(s,y,p,f,r)+\sum_{p\leq p_0,p\notin \mathcal S}j(s,y,p).
\end{align*}
There, $\gamma$ is the Euler constant, $a(x)=(3x^{-3}(\sin(x)-x\cos(x)))^2$, $\mathcal S$ is a finite (possibly empty) set of (rational) primes, $p_0$ some fixed prime and $R(\mathcal S,p_0)$,$q(s,y,p,f,r)$ and $j(s,y,p)$ are functions as defined in \cite[pp.p14-515]{ChinburgFriedman86}.
Furthermore we have the inequality
\begin{align}
\label{chfr2}
\notag \log g(k,B) &>\\
&> \log 2-\log([k_B:k])+m T(0,y)\\
\notag &\geq \log(2[U_k:U_{k}^+])-m\log 2-\log([k'_B:k])+m T(0,y).
\end{align}
\\
Moreover, if all primes in $k$ lying over $2$ that belong to $R_f(B)$ are decomposed in $K$ and functions $q(s,y,p,r,f)$ and $j(s,y,p)$ are replaced by $\hat q(s,y,p,f,r)$ and $\hat j(s,y,p)$ as in \cite[pp.514-515]{ChinburgFriedman86}, then the inequalities (\ref{chfr1}) and (\ref{chfr2}) above hold with $g(k,B)$ replaced by $\mu(k,B)$.  
\end{thm} 
\begin{remark}
The statement is slightly different from the original statement of \cite[Proposition 3.1]{ChinburgFriedman86} and is obtained by replacing Zimmert's lower bound $0.02\exp(0.46 m)$ for the regulator, which is used in \cite{ChinburgFriedman86}, by the better bound $0.0062\exp(0.738 m)$ due to Friedman, see \cite{Friedman89}
\end{remark}

Now, evaluating $T(s,y)$ for $s=1.4$ and $y=0.1$ and applying the above theorem to $K=k'_B$, $\mathcal S=\emptyset$, $p_0=31$, as in \cite{ChinburgFriedman86} but using slightly finer estimates 
$$
\frac{3-1.4}{2}\int_{0}^{\infty}\frac{1-a(\sqrt{0.1}x)}{\sinh(x)}dx<0.061530978
$$
$$
\frac{3-1.4}{4}\int_{0}^{\infty}\frac{1-a(\sqrt{0.1}x)}{\cosh^2(x/2)}dx<0.047218236
$$
computed with MAPLE, we get the inequality

\begin{equation}
g(k,B)> \psi(m,[k'_B:k]):=0.477132753\exp(0.658829093\cdot m-\frac{19.0745}{[k'_B:k]})
\end{equation}
The function $\psi(m,[k'_B:k])$ is increasing when $m$ or $[k'_B:k]$ increases, and after computing $\psi(24,1)=0.01827...$, $\psi(10,2)=0.0250\ldots$ and $\psi(3,4)=0.29244\ldots$ we conclude
\begin{lem}
Let $\Gamma\in \mathcal C(k,B)$ be an arithmetic lattice in $G_n$ such that $\chi(\Gamma)\leq 1/60=0.0166\ldots$. Then $m=[k:\QQ]\leq 23$. For $10\leq m\leq 23$ the abelian extension $k'_B$ associated with $\mathcal C(k,B)$ coincides with $k$ and for $3\leq m\leq 9$ the field $k'_B$ is at most a quadratic extension of $k$. 
\end{lem}

We now follow the steps of \cite{ChinburgFriedman86} to deal with classes $\mathcal C(k,B)$ where $m\leq 23$ and $[k'_B:k]\leq 2$ and make now use of the inequality (\ref{chfr2}) with $K=k'_B$. Again, using the estimates for $T(0,1)$ as in \cite{ChinburgFriedman86} but replacing
$$
\frac{3}{2}\int_{0}^{\infty}\frac{1-a(x)}{\sinh(x)}dx<0.71554278
$$
$$
\frac{3}{4}\int_{0}^{\infty}\frac{1-a(x)}{\cosh^2(x/2)}dx<0.492435643
$$
(we used MAPLE) we find that 

\begin{equation}
\label{mu}
\log\mu(k,B)>\xi(m,[k'_B:k]):=1.38629436-\log([k'_B:k])+0.58548009\cdot m-\frac{11.3098}{[k'_B:k]}
\end{equation}

Now, if for $\Gamma\in \mathcal C(k,B)$ we have $\chi(\Gamma)\leq 1/60$ and $10\leq m\leq 23$, then $[k'_B:k]=1$ and we compute $exp(\xi(10,1))=0.01709\ldots$ and $\exp(\xi(2,2))>1/60$. It follows
\begin{lem}
\label{mleq5}
If $\Gamma\in \mathcal C(k,B)$ satisfies $\chi(\Gamma)\leq 1/60$ then $m\leq 9$. In this case $k'_B=k$. 
\end{lem} 

In the next step we follow again \cite{ChinburgFriedman86} and evaluate the basic inequality (\ref{chfr2}) with $K=k_B$, $\mathcal S=\emptyset$, $p_0=2$ to obtain
 
\begin{equation}
\log \mu(k,B) > \rho(m,[k_B:k]):=0.693147-\log([k_B:k])+m\cdot 1.278628-\frac{11.3098}{[k_B:k]}
\end{equation}
We compute $\exp(\rho(6,1))=0.05\ldots>1/60$. This, together with the computation of $\rho(m,2)>1/60$ for $m\leq 5$ shows 
\begin{lem}
If $\mathcal C(k,B)$ contains a lattice $\Gamma$ with $\chi(\Gamma)\leq 1/60$, then $m=[k:\QQ]\leq 5$. Moreover $k_B=k$. 
\end{lem} 

\subsection{Analysis of the classes $\mathcal C(k,B)$ with $m\leq 5$}
We will now assume that $\mathcal C(k,B)$ contains a lattice $\Gamma$ with $\chi(\Gamma)\leq 1/60$. Then $m=[k:\QQ]\leq 5$ and $k_B=k$. Using the well-known lower bound $\zeta_k(2)>\zeta_{\QQ}(2m)$ and the fact that $k_B=k$ from Lemma \ref{mleq5} we know from the definition (\ref{gkb}) of $g(k,B)$ that in this case 

$$
g(k,B)>\alpha(m,d_k,t):=\frac{d_k^{3/2}\zeta_{\QQ}(2m)}{2^{2m+t-1}\pi^{2m}}.
$$
Since $t\leq m$, we observe that always $\alpha(m,d_k,t)\geq \alpha(m,d_k,m)$.

\subsubsection{Degree $m=5$} Consider first the case $m=5$. Then, the explicit computation of $\alpha(5,89000,5)=0.01732200133$ shows that the discriminant of the totally real quintic field $k$ is less then $89000$, when $g(k,B)\leq 1/60$. There are exactly $8$ such fields (see, for instance, \cite{NF} or \cite{Voighthomepage}) and for these fields we can compute the exact zeta values $\zeta_k(2)$, using PARI or SAGE.
\begin{table}[ht]
\centering
\begin{tabular}{| c | c | c | c | c | c |}
\hline
$d_k$ & defining polynomial & $\zeta_k(-1)$ & $t$ & $k_B$ & $g(k,B)$\\
\hline
$ 14641$ & $ x^5 -x^4 -4x^3 + 3x^2 +3x - 1$ & $-60/99$ & $1$ & $k$ & $0.18939\ldots$ \\
\hline
$ 24217$ & $ x^5 - 5x^3 - x^2 + 3x + 1$ & $-4/3 $ & $1$ & $k$ & $1/24$ \\
\hline
$ 36497$ & $x^5 - 2x^4 - 3x^3 + 5x^2 + x - 1$ & $-8/3$ & $1$ & $k$ & $1/12$\\
\hline
$38569 $ & $ x^5 - 5x^3 + 4x - 1$ & $-8/3$ & $1$ & $k$ & $1/12$\\
\hline
$65657$ & $x^5-x^4-5x^3+2x^2+5x+1$ & $-20/3$ & $1$ & $k$ & $0.2083\ldots$\\
\hline
$70601$ & $x^5-x^4-5x^3+2x^2+3x-1$ & $-20/3$ & $1$ & $k$ & $0.2083\ldots$\\
\hline 
$81509$ & $x^5 - x^4 - 5x^3 + 3x^2 + 5x - 2$ & $-8/3$ & $1$ & $k$ & $1/12$\\
\hline
$81589$ & $x^5-6x^3+8x-1$ & $-32/3$ & $2$ & $k$ & $1/6$\\
\hline
\end{tabular} 
\caption{$g(k,B)$-values of totally real quintic fields with discriminant $\leq 89000$ }\label{grad5}  
\end{table}

The Table \ref{grad5} shows that in degree $m=5$ all lattices in $\mathcal C(k,B)$ satisfy $\chi(\Gamma)>1/60$. 

\subsubsection{Degree $m=4$}
As in the previous case, we can compute $\alpha(m,d_k,m)$ and we obtain $\alpha(4,4800,4)=0.171830$, from which follows that we need to analyze the totally real quartic fields with discriminants $d_k<4800$. There are $16$ such fields (again consulting \cite{NF}) and the explicit computation with SAGE (see table \ref{grad4}) identifies two fields of degree $4$ which may lead to lattices of covolume less or equal $1/60$: 
\begin{itemize}
\item For $k=k_{725}$ the totally real field of discriminant $725$ we have $g(k,B)=1/120$. If $R(B)=\emptyset$, that is, $B=M_2(k_{725})$, then, as $U_{k_{725}}^+=U_{k_{725}}^2$, the Hilbert modular group $\Delta_{k_{725}}^{nu}=\PSl_2(\oo_{k_{725}})$ is maximal and we have $\chi(\PSl_2(\oo_{k_{725}}))=\mu(k,B)=1/60$. The same value is attained also by the maximal lattice $\Delta_{k_{725}}^n=N\Gamma_{\OO}^+\subset G_2$ inside $\mathcal C(k_{725},B)$, where $B$ is the division quaternion algebra ramified at two infinite places of $k_{725}$. Suppose, $B$ is ramified at exactly one infinite place. Since prime $2$ is inert in $k_{725}$, $t'=0$. Then, as the number of ramified places is even, there is one finite ramified place $v$ in $B$ contributing $(Nv-1)/2$ to $\chi(\ )$. Also, $3$ is inert in $k_{725}$ and $5=\pP_5^2$ is ramified. Analyzing the splitting behavior of small primes we conclude that a ramified finite place contributed at least $12$ to $\chi(\ )$, which by far exceeds $1/60$.
\item The other candidate is $k_{1125}$. But as $2$ is inert in $k_{1125}$, $t'=0$ and $\chi(\Gamma)\geq \mu(k,B)=2g(k,B)=1/30$ for any $\Gamma\in \mathcal C(k,B)$.  
\end{itemize}

\begin{table}[h]
\centering

\begin{tabular}{|c|c|c|c|c|c|}
\hline
$d_k$  & defining polynomial & $\zeta_k(-1)$ &  $t$ & $[k_B:k]$ & $g(k,B)$ \\
\hline
$725 $ & $x^4 - x^3 - 3x^2 + x + 1$ 	&  $2/15$ & $1$ & $1$ & $1/120$ \\
\hline
$1125$ & $x^4 - x^3 - 4x^2 + 4x + 1$	&  $4/15$ & $1$ & $1$ & $1/60$ \\
\hline
$1957$ & $x^4 - 4x^2 - x + 1$		&  $2/3$ & $1$ & $1$ & $1/24$  \\
\hline
$2000$ & $x^4 - 5x^2 + 5$ 		&  $2/3$ & $1$ & $1$ & $1/24$ \\
\hline
$2225$ & $x^4 - x^3 - 5x^2 + 2x + 4$	&  $4/5$ & $1$ & $1$ & $1/20$ \\
\hline
$2304$ & $x^4 - 4x^2 + 1$		&  $1$ & $1$ & $1$ & $1/16$  \\
\hline
$2525$ & $x^4 - 2x^3 - 4x^2 + 5x + 5$	&  $4/3$ &  $1$ & $1$ & $1/12$ \\
\hline
$2624$ & $x^4 - 2x^3 - 3x^2 + 2x + 1$	&  $1$ & $1$ & $1$ & $1/16$  \\
\hline
$2777$ & $x^4 - x^3 - 4x^2 + x + 2$	&  $4/3$ & $2$ & $1$ & $1/24$ \\
\hline
$3600$ & $x^4 - 2x^3 - 7x^2 + 8x + 1$	&  $8/5$ & $1$ & $1$ & $1/10$ \\
\hline
$3981$ & $x^4 - x^3 - 4x^2 + 2x + 1$	&  $2$ & $1$ & $1$ & $1/16$ \\
\hline
$4205$ & $x^4 - x^3 - 5x^2 - x + 1$	&  $2$ & $1$ & $1$ & $1/8$ \\
\hline
$4352$ & $x^4 - 6x^2 - 4x + 2$		& $8/3$ & $1$ & $2$ & $1/12$ \\
\hline
$4400$ & $x^4 - 7x^2 + 11$	& $34/15$ & $1$ & $2$ & $0.07\ldots$\\
\hline
$4525$ & $x^4 - x^3-7x^2+3x+9$	& $34/15$ & $1$ & $1$ & $0.14\ldots$\\
\hline
$4752$ & $x^4 - 2x^3-3x^2+4x+1$	& $8/3$ & $1$ & $2$ & $1/12$\\
\hline

\end{tabular}

\caption{$g(k,B)$-values of real quartic fields with discriminant $\leq 4800$.}\label{grad4}  
\end{table}

\subsubsection{Degree $m=3$} Here, we proceed as before: The explicit computation of $\alpha$ shows that the discriminant of $k$ is less then $300$. There are only $6$ such fields, but $g(k,B)\leq 1/60$ only for the cubic field $k_{49}$ of discriminant $49$. Suppose now $B=M_2(k_{49})$, that is $R(B)=\emptyset$. The minimal value for $\chi$ is then $\mu(k_{49},M_{2}(k_{49}))=2/168=1/84$ which is attained by the Hilbert modular group $\Delta_{k_{49}}^{nu}=\PSl_2(\oo_{k_{49}})$. Suppose now $R(B)\neq \emptyset$. If $B$ is ramified at exactly one infinite place (hence we are considering subgroups of $G_2$), then, there must be at least one finite place $v$ in $R(B)$ contributing $(Nv-1)/2$. As $2$, $3$ and $5$ are inert in $k_{49}$, $t'=0$ and this contribution is greater than $2$, hence the value of $\chi$ is greater than $1/60$. Finally, if we assume that $B$ is ramified at two infinite places, we obtain the well-known cocompact arithmetic Fuchsian group of minimal covolume.    
\begin{table}[h]
\centering

\begin{tabular}{|c|c|c|c|c|c|}
\hline
$d_k$  & defining polynomial & $\zeta_k(-1)$ &  $t$ & $[k_B:k]$ & $g(k,B)$ \\
\hline
$49$ & $x^3 - x^2 - 2x + 1$ 	&  $-1/21$ & $1$ & $1$ & $1/168$ \\
\hline
$81$ & $x^3 - 3x -1$	&  $-1/9$ & $1$ & $1$ & $1/72$ \\
\hline
$148$ & $x^3 - x^2 - 3x + 1$		&  $-1/3$ & $1$ & $1$ & $1/24$  \\
\hline
$169$ & $x^3 - x^2 -4x-1$ 		&  $-1/3$ & $1$ & $1$ & $1/24$ \\
\hline
$229$ & $x^3 - 4x -1$	&  $-2/3$ & $2$ & $1$ & $1/24$ \\
\hline
$257$ & $x^3 - x^2 -4x+3$		&  $-2/3$ & $1$ & $2$ & $1/24$  \\
\hline
\end{tabular}

\caption{$g(k,B)$-values of real quartic fields with discriminant $\leq 4800$.}\label{grad4}  
\end{table}

\subsubsection{Degree $m=2$}
In \cite{Dzambic14} we determined all the commensurability classes $\mathcal C(k,B)$ with $k$ real quadratic such that $g(k,B)\leq 1$. Analyzing the examples there we see that only possibility for $g(k,B)\leq 1/60$ is $k=\QQ(\sqrt{5})$ with $\mu(k,B)=1/60$ (note that $t'=0$ since $2$ is inert in $k$). The Hilbert modular group of the field $k$ attains this minimum. Suppose that the quaternion algebra $B$ over $k$ ramifies at some place. Then there must be at least one finite ramified place and since the primes $2$ and $3$ are inert, there must be a non-trivial contribution $(Nv-1)/2>1$ to the volume. This implies that beside the above example, there are no lattices $\Gamma\in \mathcal C(\QQ(\sqrt{5}),B)$ with $\chi(\Gamma)\leq 1/60$. 

\subsection{Minimal covolume with respect to the Euler-Poincar\'e measure}
Now, we are able to prove also the second part of the main theorem, where the lattices of minimal covolume with respect to the Euler-Poincar\'e measure are determined. Note that for an arithmetic lattice $\Gamma\subset G_n$, which belongs to a class $\mathcal C(k,B)$, we have $\mu(\Gamma)\geq 2^{n}\mu(k,B)=2^{n+t-t'}g(k,B)\geq 2^{n}g(k,B)$. Also, we have the candidates for the smallest covolumes given in the main Theorem \ref{mainthm}. Now, using the functions $\Psi(n,m,[k'_B:k])=2^n\psi(m,[k'_B:k])$, $\mathcal X(n,m,[k'_B:k])=2^n\xi(m,[k'_B:k])$, and $\mathcal R(n,m,[k'_B:k])=2^n\rho(m,[k'_B:k])$ we know that $\mu(\Gamma)>\Psi(n,m,[k':k])$, $\log \mu(\Gamma)>\mathcal X(n,m,[k'_B:k])$ and $\log \mu(\Gamma)>\mathcal R(n,m,[k'_B:k])$. Now, we only need to repeat the steps from section \ref{inequalities} to obtain the result: In the first step we compute for instance
$$
\Psi(13,13,1)=0.1\ldots>1/15=0.066\ldots
$$ 
Since $\Psi(n,m,[k'_B:k])$ is increasing, it follows that for an arithmetic lattice $\Gamma\subset G_n$ belonging to the class $\cC(k,B)$ with $\mu(\Gamma)\leq 1/15$, we have $n\leq m\leq 12$. Moreover, computing $\Psi(n,m,[k_B':k])$ for small values of $n,m,[k'_B:k]$ we know that the condition $\mu(\Gamma)\leq 1/15$ implies that for $m\geq 8$ necessarily  $[k_B':k]=1$ and for $4\leq m\leq 7$ necessarily $[k'_B:k]\leq 2$. Using the function $\mathcal X$ we find in the second step that this condition implies that $m\leq 4$ and $[k_B':k]=1$, which leads us to the list of fields treated in before. This establishes the last part of the main theorem.          
   

\bibliography{Min_Vol_Literatur}
\bibliographystyle{amsplain}


\end{document}